\documentstyle[12pt]{article}
\newtheorem{theorem}{Theorem}

\newtheorem{lemma}{Lemma}

\begin{document}

\title{\bf ON PATH INTEGRALS FOR THE HIGH-DIMENSIONAL BROWNIAN BRIDGE}
%
\author{R.~PEMANTLE\\
Department of Mathematics\\
University of Wisconsin\\
Madison, Wisconsin 53706,   U.S.A.\\
\vspace{1cm}
\and
M.D.~PENROSE\\
Department of Statistics and Applied Probability\\
University of California\\
Santa Barbara, California 93106, U.S.A.}
%
%
\maketitle
\begin{abstract}

Let $v$ be a bounded function with bounded support in $R^d$,
$d \geq 3$. Let $x,y \in R^d$. 
Let $Z(t)$ denote the path integral of $v$
along the path of a Brownian bridge in $R^d$ which runs for time
$t$, starting at $x$ and ending at $y$. As $t \to \infty$,
it is perhaps evident that
the distribution of $Z(t)$ converges weakly to that of
the sum of the integrals of $v$ along
the paths of two independent Brownian motions, starting
at $x$ and $y$ and running forever.
Here we prove
a stronger result, namely convergence of the
corresponding moment generating functions and of moments.
This result is needed for applications in physics.

\end{abstract}

\newpage

\section{Introduction}
\label{sec:I}

Suppose that $(X_s, s\geq 0)$ is a stochastic process taking values in $R^d$, and
$v$ is a Real-valued function on $R^d$. Define  the path integral 
$Z(t)$ by
\begin{equation}
Z(t) = \int_0^t v(X_s) ds.  \label{pathint}
\end{equation}
If $v$ is the indicator function of a set $Z(t)$ is the occupation 
time of that set by $(X_s)$ up to time $t$. 

If $(X_s)$ is a suitable time-homogeneous Markov process,  such as a 
Feller process,  then the distribution of $Z(t)$ can be 
studied via the Feynman-Kac
formula. See for example Williams \cite{Williams}. For example,  if $(X_s)$ is
Brownian motion and $v$ is the indicator function of the unit ball,
the distribution of $Z_{\infty}$ was obtained by Cieselski and
Taylor \cite{CT} (see also \cite{Williams}).
In the Markov process literature
the path integral $Z(t)$ is known as an {\em additive functional}.

We study here the path integral $Z(t)$ given by (\ref{pathint}),
 where the process $(X_s, 0 \leq s \leq t)$ is
a Brownian bridge, and therefore is $\em not$ time-homogeneous.
Loosely speaking, a Brownian bridge $(X_s, 0 \leq s \leq t)$ is
a Brownian motion conditioned to take some fixed value $y$ at time $t$.

Our motivation for this work came from the study of the heat equation with
added potential $v(\cdot)$; that is the equation
\begin{eqnarray}
\partial f/\partial t = [(1/2) \bigtriangledown^2 - v(x)] 
f(t;x)  \quad ( t >0, x \in R^d )
										\label{eq:[3.3]}
\end{eqnarray}
(where the Laplacian $ \bigtriangledown^2 $ acts on functions of $x$ at 
fixed $t$). This equation is closely related to
the Schr\"{o}dinger equation. In quantum Statistical Mechanics, in
$\cite{Lieb67}$ and 
$\cite{PPSP}$ for example, one studies the second virial 
coefficient for a gas consisting of particles with a two-body
interaction potential $v$, via the fundamental solution
$G(t;x,y)=G_y(t;x)$  to
the associated Bloch equation
(\ref{eq:[3.3]}); that is, the solution
$f=G_y$  to (\ref{eq:[3.3]}) with the initial condition 
   $G_y(0;x) = \delta (x-y)$.

By the Feynman-Kac formula, one can show 
$G(t;x,y)$ is equal to 
the Brownian transition density $(2 \pi t)^{-d/2}e^{-|y-x|^2/(2t)}$, 
multiplied by
the moment generating function (m.g.f.) of a random variable of the form
 $Z(t)$ given by  (\ref{pathint}) where $(X_s, 0 \leq s \leq t)$ 
is a Brownian bridge running from $x$ to $y$ in time $t$ (see \cite{PPSP},
equation (29) for details). It is this m.g.f. which
we study here.

Another motivation arises in the simulation of molecular dynamics
for a chemical reaction. See Clifford {\em et al.} \cite{CGP}.
Suppose two particles (molecules) execute independent Brownian
motions, so that one particle executes a Brownian motion
$(W_s)$ relative to the other. Suppose the particles react
(for example by mutual annihilation) at a rate given by a function 
$v$ of their separation vector $W_s$; that is, if $\tau$ denotes
  the reaction time,
\begin{equation}
P[ \tau \leq t] = E \exp \{ - \int_0^t v(W_s) dx \}. 
\end{equation}

In a simulation of the system of two particles  (maybe as part of a 
much larger system), it is natural to simulate the value of
$W_s$, say at integer values of $s$, and to interpolate between successive 
integer times $n$ and $n+1$ by taking the path of $(W_s)$ between those
times to be a Brownian bridge with endpoints given by the
simulated values of $W_n$ and $W_{n+1}$. Given that 
no reaction takes place before time $n$, and given the simulated
values of $W_n$ and $W_{n+1}$, the probability of reaction 
before time $n+1$ is again the m.g.f. of a path integral of $v$ along 
the path of a Brownian bridge.

This approach was used by Clifford and Green \cite{CG}, although
their reaction mechanism is different from the one discussed here.
The Brownian bridge could arise similarly in simulations of
an integrated cost function in Control theory.

In the physical examples discussed above, the process $(X_s)$
of (\ref{pathint}) is typically a Brownian bridge in three dimensions.
In three or more dimensions, Brownian motion is transient, and
if $v$ has bounded support, denoted $A$ say, then 
for large $t$ one might expect the
Brownian motion to spend most of its time
outside $A$, leaving  $A$ shortly after time
0 and returning shortly before time $t$ (where `shortly'
is relative to $t$). Since the time-reversed
process is also a Brownian bridge, for large t one might 
expect the distribution of $Z(t)$ to approximate to that
of 
\begin{equation}
\int_0^{\infty} v(W_x(s))ds + \int_0^{\infty} v(W^{\prime}_y(s))ds,
\end{equation}
where $W_x$ and $W^{\prime}_y$ are independent unconditioned
Brownian motions starting at $x$ and $y$ (the endpoints of the bridge).
Theorem 1 below
 makes this intuition precise; we go beyond convergence
in distribution,
proving convergence of m.g.f.s and of moments.

\section{Statement of results}

{\em Assumptions}. We assume throughout 
that $d \geq 3$ and that
$v$ is a bounded, Borel measurable function on $R^d$ with 
bounded support. Let $x(t) \in R^d$ and $y(t) \in R^d$ be defined
for $t > 0$.
For each $t>0$ let the random variable $Z(t)$ be given by the path
integral (\ref{pathint}) of $v$, along the path ($X_s,0 \leq s \leq t$)
of a $d$-dimensional Brownian bridge running
for time $t$  with $X_0=x(t)$
and $X_t=y(t)$.  We view 
$\{Z(t), t \geq 0\}$  not as a stochastic process, merely  as a 
family of probablity distributions indexed by $t$. 

For $x \in R^d$, let
 $(W_x(s),s \geq  0)$ be a Wiener process (unconstrained Brownian motion)
with $W_x(0)=x$. For $t \in [0,\infty]$, define
\begin{displaymath}
Y_x(t)= \int_0^{t} v(W_x(s))ds.
\end{displaymath}
Write $Y_x$ as an abbreviation for $Y_x(\infty)$.

Before stating our main results on 
Brownian bridge path integrals, we state some facts about 
the Brownian motion path integral $X_x$. 

\begin{lemma}
There exists $\alpha_0 >0$ such that for $|\alpha|< \alpha_0$,
\begin{equation}
E \exp \{ \alpha \int_0^{\infty} |v(W_x(s))| ds \} 
< \infty,\quad \quad  x \in R^d.
\label{L1ineq}
\end{equation}
\end{lemma}
\begin{lemma}
For every real $\alpha$, $Ee^{\alpha Y_x}$ is either finite for
all $x$ or infinite for all $x$.
In the former case, $Ee^{\alpha Y_x}$ is a locally bounded, 
continuous function of  $x$.
Also, for all nonnegative integers $k$,  $E(Y_x)^k$ is continuous in 
$x$. 
\end{lemma}
\begin{lemma}
Suppose $v$ is nonnegative and not almost everywhere zero. Then there 
exists $\alpha_1 \in (0,\infty)$ such that
\begin{eqnarray}
Ee^{\alpha Y_x}< \infty, \quad \alpha < \alpha_1, \; x \in R^d;
\label{alphaless}   \\
Ee^{\alpha Y_x}= \infty, \quad \alpha > \alpha_1, \; x \in R^d.
\label{alphamore}  
\end{eqnarray}
\end{lemma}

We can now state our main results. 
All limits below are taken as $t \rightarrow \infty$.
This limiting procedure is motivated by problems arising in 
\cite{PPSP}. Using the scaling property of Brownian motion, one can re-state
the results in terms of a limiting regime where the 
time for which the Brownian bridge runs remains fixed,
and the range (support) of $v$ shrinks. 

Convergence in distribution is implied by
 the statements  (\ref{mgfconv})
and (\ref{mgfconv1}) below on
convergence of m.g.f.s on a
neighborhood of 0 (see \cite{mgfref}).
Our results on convergence of m.g.f.s are stronger than convergence in
 distribution, and are used in  
\cite{PPSP}.

\begin{theorem}
\label{thm:ybdd}
Suppose $x(t) \rightarrow x \in R^d$, and $y(t) \rightarrow y \in R^d$.
Then 
%
%
each moment of $Z(t)$ converges to the corresponding
moment of $\; Y_x + Y^{\prime}_y$, where $Y^{\prime}_y$  is defined to 
be a random variable, independent of $Y_x$, with  the
same distribution as $Y_y$. Also, 
for 
$\alpha_0$ given by Lemma 1,  
\begin{equation}
Ee^{\alpha Z(t)} \rightarrow Ee^{\alpha (Y_x + Y_y^{\prime} )}
=Ee^{\alpha Y_x}Ee^{\alpha Y_y}, 
\quad \alpha \in (- \alpha_0, \alpha_0) .
\label{mgfconv}
\end{equation}
If $\; v$ is  nonnegative, then the convergence in
(\ref{mgfconv})
holds for every $\alpha \in (- \infty, \alpha_1)$, 
where $\alpha_1$ is given by Lemma 3. 
\end{theorem}
\begin{theorem}
\label{thm:y2liket}
Suppose $x(t) \rightarrow x \in R^d$,  while $|y(t)| \rightarrow \infty$,
and $|y(t)|^2/t$ remains bounded. Then
%
every moment of $Z(t)$ converges to the corresponding
moment of $Y_x$, and for $|\alpha|<
\alpha_0$, 
\begin{equation}
Ee^{\alpha Z(t)} \rightarrow Ee^{\alpha Y_x},
 \quad \alpha \in (-\alpha_0, \alpha_0).
\label{mgfconv1}
\end{equation}
If $v$ is nonnegative, then (\ref{mgfconv1}) holds
for every $\alpha \in (-\infty, \alpha_1)$.
\end{theorem}

To prove the Theorems,  
we need the following
extension of the continuity results of Lemma 3.
\begin{lemma}
%
(a) Suppose $(x_n)$ and $(t_n)$ are sequences in $R^d$ and $(0,\infty)$
respectively, with $x_n \to x \in R^d$ and $t_n \to \infty$ as $n \to \infty$.
Then each moment of $Y_{x_n}(t_n)$ converges to the corresponding moment
of $Y_x$. Also, 
%
%
%
\begin{equation}
E \, \exp (\alpha Y_{x_n}(t_n))
\rightarrow
E \, \exp (\alpha Y_x), \quad \alpha \in (- \alpha_0,\alpha_0).
\label{expconv}
\end{equation}
If $v \geq 0$ on $R^d$, then (\ref{expconv}) holds for all
$\alpha  \in (-\infty, \alpha_1)$. 

(b) Suppose now that $|x_n| \to \infty$. Then $Ee^{\alpha Y_{x_n}(t_n)}
\to 0$ as $n \rightarrow \infty$, for any fixed $\alpha \in 
(-\alpha_0,\alpha_0)$ and for any fixed $\alpha \in (-\infty,\alpha_1)$
if $v$ is nonnegative.
\end{lemma}
\section{Proof of Lemmas}
{\em Notation.}
 Let $q(t;x)$ $(t>0,x\in R^d)$ denote
the transition density of $d$-dimensional Brownian motion;  that is,
$q(t;x)=(2\pi t)^{-d/2} \exp(-|x|^2/(2t))$.
For ${\bf s}=(s_1,\ldots,s_k) \in (0,\infty)^k$ and
${\bf z} = (z_1,\ldots,z_k) \in (R^d)^k$, let   $q_k^x({\bf s},{\bf z})$ 
denote the joint probability 
density function at $(z_1,\ldots,z_k)$ for the positions at 
times $(s_1,\ldots,s_k)$ of a $d$-dimensional Brownian
motion starting at $x$. Let $q^{x,y}_{t,k}({\bf s},{\bf z})$
denote the corresponding joint density   of a
$d$-dimensional Brownian bridge $(X_s,0\leq s \leq t)$
 starting at $x$  and finishing at $y$.

For $t \in (0,\infty]$, write $T(k;t)$ for the `triangular' set
 $\{{\bf s}=(s_1,\ldots,s_k):0<s_1<\cdots <s_k<t\}$.
For ${\bf s} \in T(k;\infty)$,  we have
\begin{equation}
q_k^x({\bf s},{\bf z})=\prod_{j=0}^{k-1}q(\Delta_j s;\Delta_j z).
\label{eq:motiondensity}
\end{equation}
where we set $s_0=0$, $z_0=x$, 
$\Delta_j s = s_{j+1}-s_j$ and
$\Delta_j z= z_{j+1}-z_j$ ($0\leq j \leq k-1$). Similarly for the
Brownian bridge; 
for ${\bf s} \in T(k;t)$,
\begin{equation}
q_{t,k}^{x,y}({\bf s},{\bf z})=[\prod_{j=0}^k
				 q(\Delta_j s;\Delta_j z)]/q(y-x;t),
\label{eq:bridgedensity}
\end{equation}
where we set  $s_{k+1}=t$,  and 
$x_{k+1}=y$.

For ${\bf z} = (z_1,z_2,\cdots,z_k) \in (R^d)^k$, set
$v_k({\bf z})=v(z_1)v(z_2) \cdots v(z_k)$.

%
\vspace{0.15in}

{\em Proof of Lemma 1}.
We have
\begin{equation}
  E\left( \int_0^{\infty} |v(W_x(s))| \: ds \right)^k/k! 
  = \int_{T(k;\infty)}    
%
       \left( \: \int_{(R^d)^{k}} |v_k({\bf z})| q_{k}^x({\bf s},{\bf z})
d{\bf z} \right)  d{\bf s}.
     \label{momentmod}
\end{equation}
Since $v$ is bounded and has bounded support, there exists finite $K_1$
for which
\begin{equation}
       \int _0^{\infty} \int_{R^d} |v(z)| q(s;z-y) dz ds < K_1,
\quad y \in R^d.
     \label{eq:momentbound}
\end{equation}
This implies that the expression in (\ref{momentmod}) is at most
$K_1^k$; hence for  $|\alpha|<K^{-1}_1$, the
expansion of the m.g.f.  (\ref{L1ineq}) as a series
of moments is absoloutely convergent.

{\em Proof of Lemma 2}.
Set $f(x)=Ee^{\alpha Y_x}$. Suppose for some $x \in R^d$
that $f(x)< \infty$. We must show that
 $f(y)< \infty$ for all $y \in R^d$.

Let $B$ be a ball of radius $r$, centered at $x$. Take $s>0$, and
for all $y \in B$, set $\tau_y=\inf \{ t\geq 0:W_y(t)  \notin B\}$, and
$\sigma_y = \min (\tau_y,s)$.
Define the measure $\mu_y$ on $\partial B$, by 
$\mu_y (E) = P[W_y(\sigma_y) \in E]$. Let $K=\sup \{|v(z)|:z \in R^d \}$.  
By using the Markov property for $W_x$ at time $s$,
\begin{equation}
f(x) \geq  e^{-Ks} \int_B f(z) q(s;z-x) dz
\label{smpats}
\end{equation}
and by the strong Markov property at time $\sigma_x$,
\begin{equation} 
f(x) \geq e^{-Ks} \int_{\partial B} f(z) \mu_x(dz).
\label{smpatsigmax}
\end{equation}

Suppose 
$|y-x| \leq r/2$. By the
 strong Markov property at time $\sigma_y$,
\begin{equation} 
f(y) \leq e^{Ks} \left( \int_B f(z) q(s;z-y) dz +
\int_{\partial B} f(z) \mu_y(dz) \right) .
\label{smpatsigmay}
\end{equation}
Since for $s$ fixed, $q(s;z-y)/q(s;z-x)$ is uniformly bounded on $z \in B$,
and $\mu_y$ has bounded Radon-Nikodym derivative with respect to
$\mu_x$, it follows that $f(y)< \infty$, and indeed
$f$ is bounded on $\{ y :|y-x| \leq r/2 \}$. Since $r$ was arbitrary,
$f$ is everywhere finite and is locally bounded.

We now show $f$ is continuous.  In the above,
let the ball $B$ have radius 1. Let $M$ be the supremum 
of $f$ over $\partial B$, which is
finite. The second integral in (\ref{smpatsigmay}) is at most
$M P[\tau_{y} < s]$,
so if $y_n \to x$, then by  
(\ref{smpatsigmay}), 
\begin{equation} 
\limsup f(y_n) \leq e^{Ks} \left( \int_B f(z) q(s;z-x) dz +
M P[\tau_{x} < s] \right) .
\label{limsupn}
\end{equation}
Now combine (\ref{limsupn}) with (\ref{smpats}), and make $s \to 0$,
to obtain  $\limsup f(y_n) \leq f(x)$. A similar argument
gives us $\liminf f(y_n) \geq f(x)$.

To show  $E(Y_x)^k$ is continuous in $x$ for
nonnegative integer $k$, use induction on $k$.
Set
\begin{displaymath}
 f_k(x)=( Ee^{\alpha Y_x} - \sum_{j=0}^k \alpha^jE(Y_x)^j/j!) 
/ \alpha^{k+1} \, . 
\end{displaymath}
Then $f_k$ is continuous in $x$ by the inductive hypothesis,
and as $\alpha \to 0$, $f_k(x)$ converges to $E(Y_x)^{k+1}/(k+1)!$,
uniformly in $x$ (see the proof of Lemma 1).

{\em Proof of Lemma 3}. If $v$ is nonnegative, then $Ee^{\alpha Y_x}$
is clearly monotone in $\alpha$, and Lemmas 1 and 2 imply
the existence of $\alpha_1 \in (0,\infty]$ satisfying
(\ref{alphaless}) and (\ref{alphamore}).
Also, $\alpha_1$ is finite; indeed, a straightforward 
geometric series argument \cite{Durrett}
shows that for given $x$, $Ee^{\alpha Y_x}$
is infinite for large enough $\alpha$. Since we shall not use the fact 
that $\alpha_1 < \infty$ here, we omit details.

{\em Proof of Lemma 4. (a)} 
First suppose $v$ is nonnegative. 
By Lemma 2,
\begin{equation}
E \exp \{ \alpha Y_{x_n}(t_n) \} \leq E \exp \{\alpha  Y_{x_n} \}
\rightarrow E \exp \{ \alpha Y_x \}.
\label{explimsup}
\end{equation}
Also, for fixed $t_0 < \infty$, the series 
$(E|\alpha Y_{x_n}(t_0)|^k/k!,k \geq 0)$
is dominated (for all $n$) by the summable series 
$((\alpha Kt_0)^k/k!, k \geq 0)$. Each moment of $Y_{x_n}(t_0)$  converges to the corresponding
moment of $Y_x(t_0)$; to see this, write each moment as a multiple integral with 
respect to the transition density $q_k^{x_n}$ or $q_k^{x}$, 
and use a domination argument.
Hence,
\begin{equation}
E \exp \{ \alpha Y_{x_n}(t_0) \} \rightarrow
E \exp \{ \alpha Y_{x} (t_0) \},
\label{fixtimeexp}
\end{equation}
and by making $t_0 \rightarrow \infty$,
\begin{displaymath}
\liminf E \exp \{ \alpha Y_{x_n}(t_n) \} 
\geq  E \exp \{ \alpha Y_{x}\}
\end{displaymath}
so that by comparison with (\ref{explimsup})
we have (\ref{expconv}). 

Now drop the assumption that $v$ is nonnegative.
If $|\alpha| < \alpha_0$, then (\ref{expconv}) would hold
if we replaced
$v(\cdot)$ by $|v(\cdot)|$ in the definitions of 
$Y_{x_{n}}(t_n)$ and $Y_x$.
Hence by obvious term-by-term estimation on the difference
of the power series,
 $|Ee^{\alpha Y_x(t)} - Ee^{\alpha Y_x}| \to 0$
locally
uniformly in $x$, and  (\ref{expconv}) follows. 
The proof of convergence of moments
is similar.

$(b)$ This is a routine application of the strong Markov property,
using the fact that $v$ has bounded support and
$Ee^{\alpha Y_x}$ is locally bounded in $x$.

\section{Proof of Theorems}

Take an arbitrary increasing Real-valued function $u(t),t \geq 0$ with
$u(0)=0$, and with $u(t) \to \infty$ and $u(t)/t \to 0$ as $t \to \infty$.
The idea of this proof is to compare the moments of $Z(t)$
with those of $Y_{x(t)}(u(t))+Y^{\prime}_{y(t)}(u(t))$.

{\em Proof of Theorem 1}. It suffices to prove 
the results along an arbitrary sequence $(t_n)$
with $t_n \to \infty$. Write $x_n$ for $x(t_n)$, 
$y_n$ for $y(t_n)$,  and $u_n$ for $u(t_n)$. Suppose
first that $v$ is nonnegative.

Let 
$T_1(k;t_n)$ 
denote the set 
$\{ {\bf s} \in T(k;t_n): s_j \in (0,u_n)\cup (t_n-u_n,\infty),
1\leq j \leq k \}$, 
and let
 $T_2(k;t_n)$
 denote the set
$T(k;t_n) \setminus T_1(k;t_n)$.
Then
\begin{eqnarray}
E(Z(t_n))^k/k!=\int_{T(k;t_n)}\int_{(R^d)^k}
q^{x_n,y_n}_{t_n,k}({\bf s};{\bf z}) 
v_k({\bf z}) d{\bf z} d{\bf s}              
\label{zeemoment}   \\
=I_1+I_2
\label{I1I2}
\end{eqnarray}
where $I_i$ denotes the integral over $T_i(k;t_n)$ ($i=1,2$).

Now consider the sum of independent unconstrained
Brownian path integrals $Y_{x_n}(t_n)+Y^{\prime}_{y_n}(t_n)$.
 By the Binomial theorem,
\begin{displaymath}
E(Y_{x_n}(t_n)+Y^{\prime}_{y_n}(t_n))^k/k! = \sum_{j=0}^k (E
[Y_{x_n}(t_n)]^j/j!)(E[Y_{y_n}(t_n)]^{k-j}/(k-j)!)
\end{displaymath}
\begin{displaymath}
> \sum_{j=0}^k \int_{T(j;t_n)} 
\{ 
\int_{T(k-j;t_n-s_j)}
\int_{(R^d)^{k-j}} v_{k-j}({\bf z^{\prime}})
q_{k-j}^{y_n}({\bf s^{\prime}};{\bf z^{\prime}})
 d{\bf z^{\prime}} d{\bf s^{\prime}}
\}   \end{displaymath}
\begin{equation}
\times \int_{(R^d)^j} v_j({\bf z})
 q_j^{x_n} ({\bf s};{\bf z}) d{\bf z} d{\bf s}.
\label{lb:binomial}
\end{equation}
Indeed,  the expression inside the braces is less than 
$E(Y_{y_n}(t_n))^{k-j}/(k-j)!$,  since  ${\bf s^{\prime}}$ has been integrated over
$T(k-j;t_n-s_j)$ ($s_j$ is the last component of ${\bf s}$),
 rather than $T(k-j;t_n)$.

By the substitutions
\begin{displaymath}
z_{j+r}=z^{\prime}_{k+1-(j+r)},\quad s_{j+r}=t_n-s^{\prime}_{k+1-(j+r)}
\quad  (r=1,2,\ldots,k-j),
\end{displaymath}
 the right hand side of (\ref{lb:binomial}) equals
%
%
\begin{displaymath}
 \sum_{j=0}^k \int_{T(k;t_n)}    \int_{(R^d)^k} 
 v_k({\bf z}) 
 q_j^{x_n}(s_1,\ldots,s_j;z_1,\ldots,z_j)
  \times
\end{displaymath}
\begin{equation}
 \times q_{k-j}^{y_n}
(t_n-s_k,t_n-s_{k-1},\ldots,t_n-s_{j+1};z_k,z_{k-1},\ldots,z_{j+1})
 d{\bf z} d{\bf s}.
%
     \label{lb:momentsum}
\end{equation}
On inspection, we find that by (\ref{eq:bridgedensity}), 
(\ref{eq:motiondensity}) and the symmetry of $q(t;x)$ in $x$,  the integrand
in (\ref{lb:momentsum}) equals
\begin{displaymath}
v_k({\bf z}) 
q^{x_n,y_n}_{t_n,k}(s_1,\ldots,s_k;z_1,\ldots,z_k) Q_{n,j}({\bf s},{\bf z})
\end{displaymath}
where the density ratio  $Q_{n,j}({\bf s},{\bf z})$ is defined by
\begin{equation}
Q_{n,j}({\bf s},{\bf z}) =
\frac{q(t_n;y_n-x_n)}{q(\Delta_j s;\Delta_j z)} = 
\left( \frac{\Delta_j s}{t_n} \right) ^{d/2}  
\exp \left( \frac{|\Delta_j z|^2}{\Delta_j s} -
\frac{|y_n-x_n|^2}{t_n} \right). 
\label{Qdef}
\end{equation}
Hence,
\begin{equation}
E (Y_{x_n}(t_n)+Y^{\prime}_{y_n}(t_n))^k/k!
>
\sum_{j=0}^k  \int_{T(k;t_n)} \int_{(R^d)^k}
q^{x_n,y_n}_{t_n,k}({\bf s},{\bf z}) Q_{n,j}({\bf s},{\bf z})
v_k({\bf z})
d{\bf z} d{\bf s}
\label{lb:indept}
\end{equation}

A similar calculation for 
$Y_{x_n}(u_n)+Y^{\prime}_{y_n}(u_n)$ gives rise to a similar formula
to (\ref{lb:binomial}), 
with ${\bf s}^{\prime}$ now integrated
over $T(k-j;u_n)$ and $t_n$ replaced by $u_n$ elsewhere,
and with 
equality this time, provided $n$ is so big that $u_n<t_n/2$.
By use of the 
 {\em same} change of variables as before 
(that is, running time backwards from
$t_n$, not $u_n$),  we obtain
\begin{displaymath}
E(Y_{x_n}(u_n)+Y^{\prime}_{y_n}(u_n))^k/k!
=
\end{displaymath}
\begin{equation}
= \sum_{j=0}^k  \int_{T_1(k;t_n)} \int_{(R^d)^k}
{\bf 1}_j({\bf s})
q^{x_n,y_n}_{t_n,k}({\bf s},{\bf z}) Q_{n,j}({\bf s},{\bf z})
v_k({\bf z})
d{\bf z} d{\bf s},
\label{eq:indepu}
\end{equation}
where ${\bf 1}_j$ denotes the indicator function of the set
$\{ {\bf s} \in T_1(k;t_n): s_j<u_n<t_n-u_n<s_{j+1} \}$.

Take the summation inside the integral in (\ref{eq:indepu}). 
For each ${\bf s} \in T_1(k;t_n)$, the indicator function
${\bf 1}_j({\bf s})$
is 1 for exactly one value of $j$. For that value of $j$, 
and for $v_k({\bf z}) \neq 0$, 
the ratio $Q_{n,j}({\bf s},{\bf z})$ is close to 1 when 
$n$ is large,  because
of the assumptions on $x_n$, $y_n$, and $u_n$.
Indeed, we have
\begin{equation}
E(Y_{x_n}(u_n)+Y^{\prime}_{y_n}(u_n))^k/k! \sim I_1 \quad 
{\rm as}\;\;  n \to \infty,
\label{I1good}
\end{equation}
in the sense that the ratio of the two sides approaches 1, 
{\em uniformly} in $k$.

Turning to $I_2$, we have that for sufficiently large
$n$,  $\exp(-|y_n-x_n|^2/t_n)>1/2$, and hence
\begin{equation}
\sum_{j=0}^k Q_{n,j}({\bf s},{\bf z}) \geq (1/2) (k+1)^{1-d/2},
\label{convexity}
\end{equation}
since by Jensen's inequality,
$\sum_{j=0}^k(\Delta_j s/t_n)^{d/2}$  is minimised by setting $s_j=t_n/(k+1)$,
for all $j$. So for large enough $n$, 
\begin{equation}
I_2 \leq 2(k+1)^{d/2-1 }
\sum_{j=0}^k
 \int_{T_2(k;t_n)} \int_{(R^d)^k} 
  q^{x_n,y_n}_{t_n,k}({\bf s},{\bf z}) Q_{n,j}({\bf s},{\bf z})
v_k({\bf z})
d{\bf z} d{\bf s}
\end{equation}
which implies (by (\ref{eq:indepu}) and (\ref{lb:indept})) that
\begin{equation}
I_2 \leq 2(k+1)^{d/2-1}D(n,k)
\label{I_2poly}
\end{equation}
where we set 
\begin{equation}
D(n,k) = [E(Y_{x_n}(t_n)+Y^{\prime}_{y_n}(t_n))^k - 
E(Y_{x_n}(u_n)+Y^{\prime}_{y_n}(u_n))^k]/k!
\label{Ddef}
\end{equation}
By Lemma 4,  $D(n,k) \to 0$ as $n \to \infty$ ;
 combining this with
(\ref{I1good}) and (\ref{I1I2}), we obtain the desired convergence of moments.

Now consider $I_2$ as a function of $k$.
Suppose $0 \leq \alpha < \alpha_1$, and
choose $\alpha_2 \in (\alpha , \alpha_1)$. For a suitable choice of 
$k_0$, we have
\begin{equation}
\sum_{k=0}^{\infty} \alpha^k I_2(k) \leq 2 \sum_{k=0}^{k_0}\alpha^k
(k+1)^{d/2-1}(D(n,k)) 
+ \sum_{k=k_0}^{\infty}\alpha_2^k(D(n,k)).
\end{equation}
The first sum converges to zero as $n \to \infty$. 
The second is bounded above by
\begin{displaymath}
E \, \exp \{ \alpha_2(Y_{x_n}(t_n)+Y^{\prime}_{y_n}(t_n)) \}
- E \, \exp \{ \alpha_2(Y_{x_n}(u_n)+Y^{\prime}_{y_n}(u_n)) \}.
\end{displaymath}
By Lemma 4, the last expression converges to zero. This,
together with the uniform convergence in (\ref{I1good}), 
implies the desired convergence (\ref{mgfconv}) of m.g.f.s on 
$0 \leq \alpha<\alpha_1$, for $v$ nonnegative.

Now drop the assumption  $v$ is nonnegative. The proof of 
(\ref{I1good}) is unchanged. Also, (\ref{I_2poly}) would hold if
 $v$ were replaced by $|v|$
in the definitions of $I_2$ and $D(n,k)$; hence $|I_2| \to 0$, and
for
 $\alpha \in (- \alpha_0, \alpha_0)$,
$\sum_{k\geq 0} |\alpha | ^k |I_2(k)| \to 0$ as above,
giving us the desired convergence
(\ref{mgfconv}) of m.g.f.s. Finally,  this 
implies convergence in distribution, so (\ref{mgfconv}) holds for 
$\alpha < 0$ when $v$ is nonnegative.

\vspace{0.25in}

{\em Proof of Theorem 2}. We use notation from the previous proof.
First consider the case that $|y_n| \to \infty$, with $|y_n|^2/t_n \to 0$.
Then convergence in the desired sense of $Z(t_n)$ to $Y_x$ follows 
by the proof of Theorem 1 (now using Lemma 4(b)).

Now consider the case that $|y_n|^2/t_n$ is bounded away from 0
and from $\infty$. First assume $v$ is nonnegative. 
Partition the set $T(k;t_n)$ as
$T_3(k;t_n) \cup T_4(k;t_n)$, where we set
\begin{displaymath}
 T_3(k;t_n) = T(k;u_n),\quad \quad 
T_4(k;t_n) = T(k;t_n) \setminus T_3(k;t_n). 
\end{displaymath}
For $i=3,4$, let $I_i$ denote the integral in (\ref{zeemoment}), restricted
to $T_i(k;t_n)$; that is,
\begin{eqnarray}
I_i=\int_{T_i(k;t_n)}\int_{(R^d)^k}
q^{x_n,y_n}_{t_n,k}({\bf s};{\bf z}) 
v_k({\bf z}) d{\bf z} d{\bf s}.              
\label{partzeemoment}  
\end{eqnarray}
So $E(Z(t_n))^k/k!=I_3+I_4$.

A similar (easier) argument to the derivation of (\ref{eq:indepu})
gives us
\begin{equation}
E(Y_{x_n}(u_n))^k/k!
=
 \int_{T(k;u_n)} \int_{(R^d)^k}
q^{x_n,y_n}_{t_n,k}({\bf s},{\bf z}) Q_{n,k}({\bf s},{\bf z})
v_k({\bf z})
d{\bf z} d{\bf s},
\label{eq:th2guts}
\end{equation}
where $Q_{n,k}({\bf s},{\bf z})$ is given by (\ref{Qdef}), with
$j=k$ and $\Delta_k z = y_n - z_{k}$ and  $\Delta_k s = t_n - s_{k}$
as before. Since
$|y_n-z|^2/|y_n|^2 \to 1$ as $n \to \infty$, 
locally uniformly in $z$, we can write
\begin{displaymath}
Q_{n,k}({\bf s},{\bf z})= (\Delta_k s/t_n)^{d/2}
\exp \{ (|y_n|^2/t_n)[(\frac{t_n}{\Delta_k s})
g_n({\bf s},{\bf z})-h_n({\bf s},{\bf z})]
\}
\end{displaymath}
where $g_n$ and $h_n$ converge to 1,  uniformly on
${\bf s} \in T(k,t_n)$  and $\bf z$ in the support of $v_k$.

Observe that  $Q_{n,k}({\bf s},{\bf z}) \to 1$,
uniformly on ${\bf s} \in T (k;u_n)$ and ${\bf z}$ 
in the support of $v_k$. Therefore,
$I_3 \sim E(Y_{x_n}(u_n))^k$ as $n\to \infty$, uniformly in $k$.

To estimate $I_4$, note that  $|y_n|^2/t_n$ is assumed to be bounded away from
0 and $\infty$ as $n \to \infty$,  so we can find $\delta>0$ 
and $n_0 \in Z$ for which $Q_{n,k}({\bf s},{\bf z}) \geq \delta$,
for all ${\bf s} \in T(k;t_n)$, ${\bf z} \in 
{\rm Supp}(v_k)$,  $n \geq n_0$ and $k \geq 0$. Thus for
$n \geq n_0$, 
\begin{displaymath}
I_4 \leq \delta^{-1} \int_{T_4} 
q^{x_n,y_n}_{t_n,k}({\bf s};{\bf z}) Q_{n,k}({\bf s},{\bf z})
v_k({\bf z}) d{\bf z} d{\bf s}.
\end{displaymath}
A similar expression to (\ref{eq:th2guts}) holds with $u_n$
replaced by $t_n$, so 
%
\begin{displaymath}
I_4 \leq \delta^{-1} [E(Y_{x_n}(t_n))^k - E(Y_{x_n}(u_n))^k]/k!
\end{displaymath}
which implies that as $n \to \infty$,
$I_4 \to 0$ and $\sum_{k\geq 0} \alpha ^k I_4(k) \to 0$
  ($0 \leq \alpha < \alpha_1$).
 Combining the estimates obtained for $I_3$ and
$I_4$, we obtain the desired convergence results for $v$ nonnegative
and $\alpha \in [0,\alpha_1)$.
The full result follows as in the proof of
Theorem 1.


\end{document}